\documentclass[11pt]{amsart}

\usepackage{amscd,verbatim}
\usepackage{amsmath, amssymb}
\usepackage{amsfonts}
\newcommand{\de}{\partial}
\newcommand{\db}{\overline{\partial}}

\newcommand{\ddbar}{i \partial \overline{\partial}}
\newcommand{\Ric}{\mathrm{Ric}}
\newcommand{\ov}[1]{\overline{#1}}

\newcommand{\tr}[2]{\mathrm{tr}_{#1}{#2}}
\newcommand{\ti}[1]{\tilde{#1}}
\newcommand{\vp}{\varphi}
\newcommand{\ve}{\varepsilon}

\renewcommand{\leq}{\leqslant}
\renewcommand{\geq}{\geqslant}

\begin{document}
\newcounter{remark}
\newcounter{theor}
\setcounter{remark}{0}
\setcounter{theor}{1}
\newtheorem{claim}{Claim}
\newtheorem{theorem}{Theorem}[section]
\newtheorem{lemma}[theorem]{Lemma}
\newtheorem{corollary}[theorem]{Corollary}
\newtheorem{conjecture}[theorem]{Conjecture}
\newtheorem{proposition}[theorem]{Proposition}
\newtheorem{question}{question}[section]
\newtheorem{defn}{Definition}[theor]
\theoremstyle{definition}
\newtheorem{rmk}{Remark}[section]

\newenvironment{example}[1][Example]{\addtocounter{remark}{1} \begin{trivlist}
\item[\hskip
\labelsep {\bfseries #1  \thesection.\theremark}]}{\end{trivlist}}

\title[Cheng-Yau estimate for symplectic Calabi-Yau]{A Cheng-Yau type estimate for the symplectic Calabi-Yau equation}

\author{Valentino Tosatti}
\address{Courant Institute of Mathematical Sciences, New York University, 251 Mercer St, New York, NY 10012}
\email{tosatti@cims.nyu.edu}

\begin{abstract}
In the setting of Donaldson's conjecture on the Calabi-Yau equation on symplectic $4$-manifolds, we prove an a priori estimate which in the K\"ahler case resembles a classical estimate of Cheng-Yau.
\end{abstract}

\maketitle

\section{Introduction}
The Calabi Conjecture, solved by Yau in 1976 \cite{Ya}, says that given a smooth positive volume form $\sigma$ on a compact K\"ahler manifold $(M,\omega)$ with total mass $\int_M \sigma=\int_M\omega^n$, we can find a unique K\"ahler metric $\ti{\omega}$ in the same cohomology class as $\omega$ whose volume element is pointwise equal to $\sigma$, i.e $\ti{\omega}^n=\sigma$. The uniqueness statement was proved by Calabi in 1954 \cite{Ca}, and the existence follows via a continuity method from {\em a priori} $C^\infty$ estimates for a K\"ahler metric which only depend on its volume form, its cohomology class, and the underlying complex manifold.

In 2006, Donaldson \cite{Do} proposed a conjectural generalization of this result to symplectic $4$-manifolds. More precisely, suppose $(M^4,J)$ is a closed almost-complex manifold and $\omega$ is a symplectic form taming $J$, which means that
$$\omega(X,JX)>0,\quad \text{for all }X\neq 0.$$
The $(1,1)$-part of $\omega$ with respect to $J$ is denoted by $\omega^{(1,1)}$.
Let $\ti{\omega}$ be another symplectic form, cohomologous to $\omega$ and compatible with $J$, which means that it tames $J$ and furthermore
$$\ti{\omega}(X,Y)=\ti{\omega}(JX,JY).$$
 There are associated Hermitian metrics $g,\ti{g}$ which are defined by
$$g(X,Y)=\frac{1}{2}\left(\omega(X,JY)+\omega(Y,JX)\right),\quad \ti{g}(X,Y)=\ti{\omega}(X,JY)$$
with corresponding volume forms $(\omega^{(1,1)})^2$ and $\ti{\omega}^2.$
Donaldson then conjectured the following:

\begin{conjecture}\label{con}
Let $(M^4,J)$ be a closed almost-complex $4$-manifold, $\omega$ a symplectic form taming $J$, $\ti{\omega}$ a cohomologous symplectic form compatible with $J$, and $\sigma$ a smooth positive volume form with $\int_M\sigma=\int_M \omega^2$. If $\ti{\omega}$ satisfies the  Calabi-Yau equation
\begin{equation}\label{scy}
\ti{\omega}^2=\sigma,
\end{equation}
then for any $k\geq 0$, we can bound $\|\ti{\omega}\|_{C^k(M,g)}$ by a constant that depends only on $k$ and on bounds on $\sigma,\omega, (M,J)$.
\end{conjecture}
When $J$ is integrable and $\omega$ is compatible with $J$, we are on a K\"ahler surface and Conjecture \ref{con} follows from the aforementioned theorem of Yau \cite{Ya}. If solved, Donaldson's conjecture would have striking applications in symplectic geometry, see \cite{Do, TW2}. Donaldson's conjecture was first investigated by Weinkove \cite{We} in the case when $\omega$ is also compatible with $J$, where he showed that it holds provided $J$ is close to being integrable. In \cite{TWY}, Tosatti-Weinkove-Yau showed that Conjecture \ref{con} in general would follow if one could prove a bound
\begin{equation}\label{c0}
\tr{g}{\ti{g}}\leq C,
\end{equation}
where $C$ depends only on $\sigma,\omega,(M,J)$ (we will call such constants uniform), and that \eqref{c0} can indeed be established when $g$ has nonnegative curvature in a suitable sense. Further progress on Donaldson's conjecture has since proceeded in two main directions: proving it on explicit examples \cite{TW3, FLSV, BFV,BFV2,TW,Ve}, and for general manifolds reducing the bound \eqref{c0} to bounding an ``almost-K\"ahler potential'' function. For this, following \cite{We,TWY}, one defines a function $\vp\in C^\infty(M,\mathbb{R})$ by
\begin{equation}\label{lapz}
\Delta_{\ti{g}}\vp=2-\tr{\ti{g}}{g}, \quad \sup_M\vp=0,
\end{equation}
which is uniquely determined and in the K\"ahler case would satisfy $\ti{\omega}=\omega+\ddbar\vp$, i.e. $\vp$ would be a familiar K\"ahler potential. Estimate \eqref{c0} was then successively reduced to proving uniform bounds for $\int_M e^{-\alpha\vp}\omega^2$ (for some $\alpha>0$) in \cite{TWY} in 2007, for $\int_M|\vp|\omega^2$ in \cite[Rmk 3.1]{TW4} in 2009 (with a very recent new proof in \cite{GP}), for {\em any} integral bound for $\vp$ in \cite{TW} in 2016, and lastly to a uniform positive lower bound for the Lebesgue measure of $\{\vp\geq -C\}$ for any given uniform $C$ in \cite{TW}.

Nevertheless, none of these results use the crucial assumption that $\ti{\omega}$ and $\omega$ are cohomologous: if this is not the case then one has to replace the constant $2$ in \eqref{lapz} with another suitable constant, namely
\begin{equation}\label{lapz2}
\Delta_{\ti{g}}\vp=\frac{2\int_M\omega\wedge\ti{\omega}}{\int_M\ti{\omega}^2}-\tr{\ti{g}}{g}, \quad \sup_M\vp=0,
\end{equation}
but there is no essential difference in all the above results, which still show that to obtain uniform $C^\infty$ bounds for $\ti{\omega}$ (and a uniform lower bound $\ti{\omega}\geq C^{-1}\omega$) it suffices to have for example any uniform integral bound for $\vp$ (or just a uniform lower bound for the measure of a superlevel set), provided the constant $\frac{2\int_M\omega\wedge\ti{\omega}}{\int_M\ti{\omega}^2}$ is uniformly bounded above. Consider the K\"ahler case, where $J$ is integrable and $\omega$ is compatible with it, and degenerate $[\ti{\omega}]$ to a limiting class $[\alpha]$ on the boundary of the K\"ahler cone with $\int_M\alpha^2>0$ (there are many such examples, see e.g. \cite{To}). In this case the constant in \eqref{lapz2} remains uniformly bounded, but the corresponding solution metrics $\ti{\omega}$ cannot converge smoothly since $[\alpha]$ does not contain any K\"ahler metric, hence in this case no uniform integral/measure bound for $\vp$ solving \eqref{lapz2} can hold. In other words, if one wants to use the almost-K\"ahler potential $\vp$ to prove \eqref{c0}, then some new idea has to be used which makes crucial use of the cohomological assumption.\\

In this note we take a different approach. Because of the cohomological assumption, we can write
\begin{equation}\label{eq}
\ti{\omega}=\omega+da,
\end{equation}
for some real $1$-form $a$. We are free to modify $a$ by adding $df$ for any smooth function $f$, and if we choose $f$ solving the elliptic equation $d^*_{\ti{g}}df=-d^*_{\ti{g}}a$, we can then assume without loss that
\begin{equation}\label{gauge}
d^*_{\ti{g}}a=0.
\end{equation}
The $1$-form $a$ satisfying \eqref{eq} and \eqref{gauge} is then uniquely determined modulo the addition of a $\ti{g}$-harmonic $1$-form, and we can fix this ambiguity for example by requiring that $a$ be $L^2(\ti{g})$-orthogonal to the space of such forms. The Calabi-Yau equation \eqref{scy} is then a nonlinear elliptic system for the $1$-form $a$ (when $\dim M=4$). Our main result is then the following:
\begin{theorem}\label{main}
In the above setting, we have the estimate
\begin{equation}\label{desi}
\tr{g}{\ti{g}}\leq C(1+\sup_M |a|^2_{\ti{g}}),
\end{equation}
where $C$ is a uniform constant.
\end{theorem}
We make a few remarks about this result.\\

\noindent 1. It is important for our arguments that the gauge-fixing condition \eqref{gauge} and the norm on the RHS of \eqref{desi} are both with respect to $\ti{g}$ as opposed to $g$. Changing \eqref{gauge} to $d^*_ga=0$ does not seem completely out of the question, but replacing $|a|^2_{\ti{g}}$ with $|a|^2_g$ in \eqref{desi} does not seem feasible.\\

\noindent 2. In the K\"ahler case, i.e. when $J$ is integrable and $\omega$ is compatible with $J$, we actually have that $a=d^c\vp$ where $\vp$ is the usual K\"ahler potential so that $\ti{\omega}=\omega+\ddbar\vp$. Indeed $dd^c\vp=\ddbar\vp$, and
$$d^*_{\ti{g}}d^c\vp=-*_{\ti{g}}d*_{\ti{g}}d^c\vp=*_{\ti{g}}d(d\vp\wedge\ti{\omega})=0,$$
and for any $\ti{g}$-harmonic $1$-form $\alpha$,
$$\int_M d^c\vp\wedge *_{\ti{g}}\alpha=\int_M d^c\vp \wedge J(\alpha)\wedge\ti{\omega}=\int_M d\vp\wedge\alpha\wedge\ti{\omega}=0.$$
Thus, in the K\"ahler case, estimate \eqref{desi} becomes
\begin{equation}\label{desi2}
\tr{g}{\ti{g}}\leq C(1+\sup_M |\nabla\vp|^2_{\ti{g}}),
\end{equation}
which in the local setting of pseudoconvex domains in $\mathbb{C}^n$ was proved by Cheng-Yau \cite[Proposition 7.1]{CY}. Tracing through our arguments easily shows that in the K\"ahler case \eqref{desi2} holds in all dimensions, and it does not use any knowledge about $L^\infty$ a priori bounds for $\vp$ (which were famously proved by Yau \cite{Ya} using Moser iteration).\\

\noindent 3. An estimate of the different but related form
\begin{equation}\label{desi3}
\tr{g}{\ti{g}}\leq C(1+\sup_M |\nabla\vp|^2_{g}),
\end{equation}
is known for different geometric PDEs on K\"ahler manifolds, including for example the complex Hessian equation \cite{HMW} and the Monge-Amp\`ere equation for $(n-1)$-psh functions \cite{TW5}, see also \cite{Sz}. In these settings, \eqref{desi3} is then used together with a blowup argument and a Liouville theorem \cite{DK, TW5, Sz} to show that $\sup_M |\nabla\vp|^2_{g}\leq C$ and hence $\tr{g}{\ti{g}}\leq C.$ However, the fact that $\ti{g}$ appears on the RHS of \eqref{desi} instead of $g$ makes it ill-suited to blowup arguments.\\

The proof of Theorem \ref{main} is done by a maximum principle argument, and interestingly it uses crucially that the dimension of $M$ is $4$ (while the aforementioned results in \cite{TWY,TW} apply in all even dimensions $2n\geq 4$), except in the K\"ahler case where our argument works in all dimensions.\\

\noindent {\bf Acknowledgments. }The author owes many thanks to Guido De Philippis for many related discussions. The author was partially supported by NSF grant DMS-2404599.

\section{Proof of Theorem \ref{main}}
We work in the setting described in the Introduction. Define a smooth function $F$ on $M$ by
$$\sigma=e^{F}(\omega^{(1,1)})^2,$$
so that in any local chart the Calabi-Yau equation \eqref{scy} can be written as
\begin{equation}\label{det}
\det(\ti{g})=e^F\det(g).
\end{equation}
Since $g,\ti{g}$ are both Hermitian with respect to $J$, and $\dim M=4$, it follows from \eqref{det} that we have
\begin{equation}\label{2d}
\tr{\ti{g}}{g}=e^{-F}\tr{g}{\ti{g}},
\end{equation}
which we will use repeatedly, often without mention.

As in \cite{TWY}, we will use covariant derivatives with respect the Chern connection $\nabla$ of $g$ (also known as ``canonical connection'', see \cite[\S 2]{TWY}). We recall two estimates proved in \cite{TWY}. The first one is the differential inequality from \cite{TWY} (which can be extracted from the proof of Lemma 3.2 there)
\begin{equation}\label{ein}
\Delta_{\ti{g}}\tr{g}{\ti{g}}\geq |\nabla\ti{g}|^2_{g,\ti{g}}-C\tr{g}{\ti{g}}\,\tr{\ti{g}}{g}-C\geq |\nabla\ti{g}|^2_{g,\ti{g}}-C(\tr{g}{\ti{g}})^2-C,
\end{equation}
using \eqref{2d}, where the norm $|\cdot |^2_{g,\ti{g}}$ uses both $g$ and $\ti{g}$, see \eqref{mixed} below for the definition. The second one is the following Cauchy-Schwarz type inequality \cite[(3.20)]{TWY}
\begin{equation}\label{ein1}
|\nabla\ti{g}|^2_{g,\ti{g}}\geq \frac{|\nabla\tr{g}{\ti{g}}|^2_{\ti{g}}}{\tr{g}{\ti{g}}}.
\end{equation}
Given these preliminaries, the main claim is then the following:
\begin{proposition}\label{prop}
There is a uniform $C>0$ such that for any small $\ve>0$ we have the differential inequality
\begin{equation}\label{zwei}
\Delta_{\ti{g}}|a|^2_{\ti{g}}\geq \frac{1}{2}|\ti{\nabla}a|^2_{\ti{g}}+\frac{1}{2}|\ov{\ti{\nabla}}a|^2_{\ti{g}}+
\frac{1}{C}(\tr{g}{\ti{g}})^2-\ve|\nabla\ti{g}|^2_{g,\ti{g}}-\frac{C}{\ve^2}|a|^4_{\ti{g}}-C.
\end{equation}
\end{proposition}
In fact, as follows from the proof, the constants $\frac{1}{2}$ on the RHS of \eqref{zwei} can be taken to be $1-\delta$ for any given $\delta>0$, at the expense of making $C$ larger. \\

\noindent {\bf Proof of Theorem \ref{main} assuming Proposition \ref{prop}.}
Let us first assume that Proposition \ref{prop} holds, and use it to prove \eqref{desi}. For this,
combining \eqref{ein} and \eqref{zwei}, and throwing away two positive terms, we have
$$\Delta_{\ti{g}}(\tr{g}{\ti{g}}+A|a|^2_{\ti{g}})\geq|\nabla\ti{g}|^2_{g,\ti{g}}-C_0(\tr{g}{\ti{g}})^2-C+\frac{A}{C_1}(\tr{g}{\ti{g}})^2
-\ve A|\nabla\ti{g}|^2_{g,\ti{g}}-\frac{CA}{\ve^2}|a|^4_{\ti{g}}-CA,$$
where $C_0,C_1$ are uniform constants, and now if we choose $A=C_1(C_0+1)$ and $\ve=\frac{1}{A}$, we obtain
$$\Delta_{\ti{g}}(\tr{g}{\ti{g}}+A|a|^2_{\ti{g}})\geq(\tr{g}{\ti{g}})^2
-C|a|^4_{\ti{g}}-C,$$
and the maximum principle concludes the proof of \eqref{desi}.\\

\noindent {\bf Beginning of proof of Proposition \ref{prop}.}
We now proceed to prove Proposition \ref{prop}. As in \cite{TWY} we also use the Chern connection $\ti{\nabla}$ of $\ti{g}$, and the formalism of moving frames. Thus, we work in a local $\ti{g}$-unitary frame $\ti{\theta}^i,i=1,2,$ which are $(1,0)$-forms, and we have the first structure equations
$$d\ti{\theta}^i=-\ti{\theta}^i_j\wedge\ti{\theta}^j+\ti{\Theta}^i,$$
where $\ti{\theta}^i_j$ are the connection $1$-forms and $\ti{\Theta}^i$ are the torsion forms, which are of type $(0,2)$ and equal to
$$\ti{\Theta}^i=\ti{N}^i_{\ov{j}\ov{k}}\ov{\ti{\theta}^j}\wedge\ov{\ti{\theta}^k},$$
where $\ti{N}^i_{\ov{j}\ov{k}}$ are the components of the Nijenhuis tensor of $J$ (and are skew-symmetric in $j,k$). The second structure equations read
$$d\ti{\theta}^i_j=-\ti{\theta}^i_k\wedge\ti{\theta}^k_j+\ti{\Omega}^i_j,$$
where $\ti{\Omega}^i_j$ are the curvature $2$-forms. They can be written as
$$\ti{\Omega}^i_j=\ti{R}^i_{jk\ov{\ell}}\ti{\theta}^k\wedge\ov{\ti{\theta}^\ell}+\ti{K}^i_{jk\ell}\ti{\theta}^k\wedge\ti{\theta}^\ell+\ti{K}^i_{j\ov{k}\ov{\ell}}\ov{\ti{\theta}^k}\wedge\ov{\ti{\theta}^\ell},$$
and the Ricci curvature form is defined to be
$$\Ric(\ti{g})=\frac{\sqrt{-1}}{2\pi}\ti{\Omega}^i_i,$$
where here and in the following we will mostly omit summation signs. It is well-known that $\Ric(\ti{g})$ is a closed real $2$-form which represents the first Chern class $c_1(M,J)$ of the complex vector bundle $(T^\mathbb{R}M,J)$. Its $(1,1)$-part is given by
$$\frac{\sqrt{-1}}{2\pi}\ti{R}_{k\ov{\ell}}\ti{\theta}^k\wedge\ov{\ti{\theta}^\ell},\quad \ti{R}_{k\ov{\ell}}=\ti{R}^i_{ik\ov{\ell}}.$$
Its relation with the Calabi-Yau equation is this \cite[(3.16)]{TWY}: if $g,\hat{g}$ are two $J$-Hermitian metrics then we have
\begin{equation}\label{ddc}
\frac{1}{2\pi}dd^c\log\frac{\det(g)}{\det(\hat{g})}=\Ric(\hat{g})-\Ric(g),
\end{equation}
where for a real-valued function $f$ we define $(1,0)$ and $(0,1)$ forms $\de f,\db f$ by $df=\de f + \db f$, and define $d^c f=\frac{\sqrt{-1}}{2}(\db f - \de f)$.
A different sign convention was used in \cite{TWY}, so the term $dd^c f$ here is equal to the term $-\frac{1}{2}d(Jdf)$ in \cite{TWY}. The trace of the $(1,1)$-part of $dd^c f$ with respect to $\ti{g}$ equals its Laplacian $\Delta_{\ti{g}} f$, see \cite[Lemma 2.5]{TWY}, which is equal to the usual Laplace-Beltrami operator of $\ti{g}$ (up to a factor of $\frac{1}{2}$). It can be written as
$$\Delta_{\ti{g}}f=2\frac{\ti{\omega}\wedge dd^c f}{\ti{\omega}^2}.$$
If we have a tensor $T$ then we can consider its covariant derivatives $\ti{\nabla} T$ and $\ov{\ti{\nabla}}T$. For example, if in our frame $T$ has components $\ti{T}^i_{j\ov{k}}$ (as an example) then its covariant derivatives $\ti{\nabla} T$ and $\ov{\ti{\nabla}}T$ have components $\ti{T}^i_{j\ov{k},p}$ and $\ti{T}^i_{j\ov{k},\ov{p}}$ respectively, which can be obtained as follows:
$$\ti{T}^i_{j\ov{k},p}\ti{\theta}^p+\ti{T}^i_{j\ov{k},\ov{p}}\ov{\ti{\theta}^p}=d\ti{T}^i_{j\ov{k}}+\ti{T}^q_{j\ov{k}}\ti{\theta}^i_q-\ti{T}^i_{q\ov{k}}\ti{\theta}^q_j-\ti{T}^i_{j\ov{q}}\ov{\ti{\theta}^q_k},$$
and similarly for tensors of other types.
We can easily compute \cite[(2.28)]{TWY} that if $f$ is a real-valued function then we have
\begin{equation}\label{glz}
dd^c f = \ti{f}_{\ov{i}}\ov{\ti{\Theta}^i}+\ti{f}_{i\ov{j}}\ti{\theta}^i\wedge\ov{\ti{\theta}^j}+\ti{f}_i\ti{\Theta}^i.
\end{equation}
To simplify computations, denote by $\alpha=a^{(1,0)}$, so that
$$a=\alpha+\ov{\alpha},$$
and we have
$$|a|^2_{\ti{g}}=2|\alpha|^2_{\ti{g}},$$
where we are using $\ti{g}$ here as a Riemannian metric on $T^{\mathbb{R}}M$ on the LHS and as a Hermitian metric on $T^{1,0}M$ on the RHS, and the equality follows from the fact that the Riemannian metric is $J$-invariant.
In our frame we can write  $\alpha=\ti{\alpha}_i\ti{\theta}^i$, so that
$$|\alpha|^2_{\ti{g}}=|\ti{\alpha}_i|^2.$$
We then compute
$$d|\alpha|^2_{\ti{g}}=\ti{\alpha}_{i,j}\ov{\ti{\alpha}_i}\ti{\theta}^j+\ti{\alpha}_{i,\ov{j}}\ov{\ti{\alpha}_i}\ov{\ti{\theta}^j}+\ti{\alpha}_i\ov{\ti{\alpha}_{i,\ov{j}}}\ti{\theta}^j+\ti{\alpha}_i\ov{\ti{\alpha}_{i,j}}\ov{\ti{\theta}^j},$$
hence
$$\de|\alpha|^2_{\ti{g}}=\ti{\alpha}_{i,j}\ov{\ti{\alpha}_i}\ti{\theta}^j+\ti{\alpha}_i\ov{\ti{\alpha}_{i,\ov{j}}}\ti{\theta}^j,$$
$$(d\de|\alpha|^2_{\ti{g}})^{(1,1)}=(\ti{\alpha}_{i,j\ov{k}}\ov{\ti{\alpha}_i}+\ti{\alpha}_{i,j}\ov{\ti{\alpha}_{i,k}}+\ti{\alpha}_{i,\ov{k}}\ov{\ti{\alpha}_{i,\ov{j}}}+\ti{\alpha}_i\ov{\ti{\alpha}_{i,\ov{j}k}})\ov{\ti{\theta}^k}\wedge\ti{\theta}^j$$
and since $dd^c |\alpha|^2_{\ti{g}}=-\sqrt{-1}d\de|\alpha|^2_{\ti{g}}$, this implies that
\begin{equation}\label{lapl}
\Delta_{\ti{g}} |\alpha|^2_{\ti{g}}=|\ti{\alpha}_{i,j}|^2+|\ti{\alpha}_{i,\ov{j}}|^2+\ti{\alpha}_{i,k\ov{k}}\ov{\ti{\alpha}_i}+\ti{\alpha}_i\ov{\ti{\alpha}_{i,\ov{k}k}}.
\end{equation}
Next, we need a commutation relation for covariant derivatives of $\alpha$. We start with the definition
$$d\ti{\alpha}_i=\ti{\alpha}_j\ti{\theta}^j_i+\ti{\alpha}_{i,j}\ti{\theta}^j+\ti{\alpha}_{i,\ov{j}}\ov{\ti{\theta}^j},$$
and applying $d$ again
\begin{equation}\label{null}
0=dd\ti{\alpha}_i=\ti{\alpha}_j\ti{\Omega}^j_i+\ti{\alpha}_{i,jp}\ti{\theta}^p\wedge\ti{\theta}^j+\ti{\alpha}_{i,j\ov{p}}\ov{\ti{\theta}^p}\wedge\ti{\theta}^j+\ti{\alpha}_{i,j}\ti{\Theta}^j+
\ti{\alpha}_{i,\ov{j}p}\ti{\theta}^p\wedge\ov{\ti{\theta}^j}+\ti{\alpha}_{i,\ov{j}\ov{p}}\ov{\ti{\theta}^p}\wedge\ov{\ti{\theta}^j}+\ti{\alpha}_{i,\ov{j}}\ov{\ti{\Theta}^j},
\end{equation}
and the $(1,1)$-part of \eqref{null} says
$$\ti{\alpha}_j\ti{R}^j_{ik\ov{\ell}}\ti{\theta}^k\wedge\ov{\ti{\theta}^\ell}-\ti{\alpha}_{i,k\ov{\ell}}\ti{\theta}^k\wedge\ov{\ti{\theta}^\ell}+\ti{\alpha}_{i,\ov{\ell}k}\ti{\theta}^k\wedge\ov{\ti{\theta}^\ell}=0,$$
and so we obtain the commutation relation
\begin{equation}\label{comm1}
\ti{\alpha}_{i,k\ov{\ell}}=\ti{\alpha}_{i,\ov{\ell}k}+\ti{\alpha}_j\ti{R}^j_{ik\ov{\ell}},
\end{equation}
exactly like in the K\"ahler case. Taking the $(0,2)$-part of \eqref{null} gives
$$\ti{\alpha}_j\ti{K}^j_{i\ov{k}\ov{\ell}}\ov{\ti{\theta}^k}\wedge\ov{\ti{\theta}^\ell}+\ti{\alpha}_{i,j}\ti{N}^j_{\ov{k}\ov{\ell}}\ov{\ti{\theta}^k}\wedge\ov{\ti{\theta}^\ell}
+\ti{\alpha}_{i,\ov{\ell}\ov{k}}\ov{\ti{\theta}^k}\wedge\ov{\ti{\theta}^\ell}=0,$$
which, after skew-symmetrizing $\ti{\alpha}_{i,\ov{\ell}\ov{k}}$ in $k,\ell$, gives the commutation relation
\begin{equation}\label{comm2}
\ti{\alpha}_{i,\ov{\ell}\ov{k}}=\ti{\alpha}_{i,\ov{k}\ov{\ell}}-2\ti{\alpha}_j\ti{K}^j_{i\ov{k}\ov{\ell}}-2\ti{\alpha}_{i,j}\ti{N}^j_{\ov{k}\ov{\ell}}.
\end{equation}
Using \eqref{comm1} in \eqref{lapl} gives
$$\Delta_{\ti{g}} |\alpha|^2_{\ti{g}}=|\ti{\alpha}_{i,j}|^2+|\ti{\alpha}_{i,\ov{j}}|^2+2\mathrm{Re}(\ti{\alpha}_{i,\ov{k}k}\ov{\ti{\alpha}_i})+\ti{\alpha}_j\ov{\ti{\alpha}_i}\ti{R}^j_{ik\ov{k}},$$
and recalling from \cite[(2.21)]{TWY} that
\begin{equation}\label{comm}
\ti{R}^j_{ik\ov{k}}=\ti{R}_{i\ov{j}}-4\ti{N}^q_{\ov{p}\ov{j}}\ov{\ti{N}^{p}_{\ov{q}\ov{i}}}-4\ti{N}^{p}_{\ov{q}\ov{j}}\ov{\ti{N}^{i}_{\ov{p}\ov{q}}},
\end{equation}
we obtain
$$\Delta_{\ti{g}} |\alpha|^2_{\ti{g}}=|\ti{\alpha}_{i,j}|^2+|\ti{\alpha}_{i,\ov{j}}|^2+2\mathrm{Re}(\ti{\alpha}_{i,\ov{k}k}\ov{\ti{\alpha}_i})+\ti{R}_{i\ov{j}}\ti{\alpha}_j\ov{\ti{\alpha}_i}-4\ti{\alpha}_j\ov{\ti{\alpha}_i}\ti{N}^q_{\ov{p}\ov{j}}\ov{\ti{N}^{p}_{\ov{q}\ov{i}}}
-4\ti{\alpha}_j\ov{\ti{\alpha}_i}\ti{N}^{p}_{\ov{q}\ov{j}}\ov{\ti{N}^{i}_{\ov{p}\ov{q}}}
.$$
To deal with the term $\ti{R}_{i\ov{j}}\ti{\alpha}_j\ov{\ti{\alpha}_i}$, observe that differentiating the PDE \eqref{det} and using \eqref{ddc} we have that $\Ric(\ti{g})=\Ric(g)-\frac{1}{2\pi}dd^c F,$ which is a fixed background tensor (independent of $\ti{g}$), hence we can write $\ti{R}_{i\ov{j}}=\ti{T}_{i\ov{j}}$ where $T=\sqrt{-1}\ti{T}_{i\ov{j}}\ti{\theta}^i\wedge\ov{\ti{\theta}^j}$ is some fixed tensor. This implies that
$$\ti{R}_{i\ov{j}}\ti{\alpha}_j\ov{\ti{\alpha}_i}\geq -|T|_{\ti{g}}|\alpha|^2_{\ti{g}}\geq -C\tr{g}{\ti{g}}|\alpha|^2_{\ti{g}}.$$
Likewise, if we let now $\ti{T}_{i\ov{j}}=\ti{N}^q_{\ov{p}\ov{j}}\ov{\ti{N}^{p}_{\ov{q}\ov{i}}}$, then $T=\sqrt{-1}\ti{T}_{i\ov{j}}\ti{\theta}^i\wedge\ov{\ti{\theta}^j}$ is some fixed tensor (independent of $\ti{g}$), and so by the same logic we can bound
$$-4\ti{\alpha}_j\ov{\ti{\alpha}_i}\ti{N}^q_{\ov{p}\ov{j}}\ov{\ti{N}^{p}_{\ov{q}\ov{i}}}\geq -C\tr{g}{\ti{g}}|\alpha|^2_{\ti{g}}.$$
Unfortunately this doesn't work directly for the last term
\begin{equation}\label{list}
-4\ti{\alpha}_j\ov{\ti{\alpha}_i}\ti{N}^{p}_{\ov{q}\ov{j}}\ov{\ti{N}^{i}_{\ov{p}\ov{q}}},
\end{equation}
since the corresponding tensor $T=\sqrt{-1}\ti{N}^{p}_{\ov{q}\ov{j}}\ov{\ti{N}^{i}_{\ov{p}\ov{q}}} \ti{\theta}^i\wedge\ov{\ti{\theta}^j}$ depends also on $\ti{g}$. To understand more precisely how to bound this, denote by $\theta^i$ a local unitary frame for the background metric $g$, so that we can write
$$\ti{\theta}^i=a^i_j\theta^j,\quad \theta^i=b^i_j\ti{\theta}^j,$$
where the local matrices $(a^i_j)$ and $(b^i_j)$ are inverses of each other, i.e. $a^i_j b^k_i=\delta_{jk}$. We can express the components of the Nijenhuis tensor with respect to the frame $\theta^i$ as
$$N^i_{\ov{j}\ov{k}}=\ti{N}^p_{\ov{q}\ov{r}}\ov{a^q_j}\ov{a^r_k}b^i_p,\quad\text{hence}\ \ti{N}^i_{\ov{j}\ov{k}}=N^p_{\ov{q}\ov{r}}\ov{b^q_j}\ov{b^r_k}a^i_p,$$
and the components $N^i_{\ov{j}\ov{k}}$ are all uniformly bounded.
Similarly the components of $(1,0)$-form $\alpha$ are given by $\alpha_i=a^j_i \ti{\alpha}_j$, and so we can write
$$4\ti{\alpha}_j\ov{\ti{\alpha}_i}\ti{N}^{p}_{\ov{q}\ov{j}}\ov{\ti{N}^{i}_{\ov{p}\ov{q}}}=4\alpha_w\ov{\alpha_h}b^w_j\ov{b^h_i}N^{k}_{\ov{\ell}\ov{r}}\ov{N^{t}_{\ov{u}\ov{v}}} \ov{b^\ell_q}\ov{b^r_j}a^p_k b^u_p b^v_q\ov{a^i_t}=4\alpha_w\ov{\alpha_h}b^w_jN^{k}_{\ov{\ell}\ov{r}}\ov{N^{h}_{\ov{k}\ov{v}}} \ov{b^\ell_q}\ov{b^r_j} b^v_q.$$
and working at an arbitrary point we may choose our unitary frames so that at this point we have
$$a^i_j=\sqrt{\lambda_i}\delta_{ij},$$
where $\lambda_1,\lambda_2>0$ are the eigenvalues of the Hermitian metric $\ti{g}$ with respect to $g$. This implies that $$b^i_j=\frac{1}{\sqrt{\lambda_j}}\delta_{ij},$$ and so at our point our term \eqref{list} simplifies to
$$-4\sum_{j,q=1}^2\lambda_j^{-1}\lambda_q^{-1}\alpha_j\ov{\alpha_h}N^{k}_{\ov{q}\ov{j}}\ov{N^{h}_{\ov{k}\ov{q}}},$$
and since $N^{k}_{\ov{q}\ov{j}}$ is skew-symmetric in $q,j$, only the terms in the sum with $j\neq q$ survive, and so this equals
$$-4\lambda_1^{-1}\lambda_2^{-1}\alpha_1\ov{\alpha_h}N^{k}_{\ov{2}\ov{1}}\ov{N^{h}_{\ov{k}\ov{2}}}-4\lambda_2^{-1}\lambda_1^{-1}\alpha_2\ov{\alpha_h}N^{k}_{\ov{1}\ov{2}}\ov{N^{h}_{\ov{k}\ov{1}}},$$
but from the Calabi-Yau equation \eqref{scy} we have that
\begin{equation}\label{cy}
\lambda_1\lambda_2=e^F,
\end{equation}
which is uniformly bounded, and so our term \eqref{list} can be bounded below by
$$-C\sum_h |\alpha_h|^2=-C|\alpha|^2_g\geq -C\tr{g}{\ti{g}}|\alpha|^2_{\ti{g}},$$
and combining all of the above gives
$$\Delta_{\ti{g}} |\alpha|^2_{\ti{g}}\geq |\ti{\alpha}_{i,j}|^2+|\ti{\alpha}_{i,\ov{j}}|^2+2\mathrm{Re}(\ti{\alpha}_{i,\ov{k}k}\ov{\ti{\alpha}_i})-C\tr{g}{\ti{g}}|\alpha|^2_{\ti{g}}.$$
Observe that here we have used crucially that $\dim M=4$ so that there are only two eigenvalues. This fact will be also used further below.\\

\noindent {\bf The main claim \eqref{vier}.}
If we denote by $\Delta^H_{\ti{g}} a=-(dd^*_{\ti{g}}a+d^*_{\ti{g}}da)$ the $\ti{g}$-Hodge Laplacian of $a=\alpha+\ov{\alpha}$, then the main claim is that
\begin{equation}\label{vier}
2\mathrm{Re}(\ti{\alpha}_{i,\ov{k}k}\ov{\ti{\alpha}_i})\geq \frac{1}{2}\ti{g}(\Delta^H_{\ti{g}} a,a)-\frac{1}{2}|\ti{\alpha}_{i,j}|^2-C\tr{g}{\ti{g}}|\alpha|^2_{\ti{g}}-C|\nabla\ti{g}|_{g,\ti{g}}|\alpha|^2_{\ti{g}}.
\end{equation}

\noindent {\bf End of proof of Proposition \ref{prop} assuming \eqref{vier}.}
Assuming that \eqref{vier} holds, let us complete the proof of Proposition \ref{prop}. Plugging in, we get
$$\Delta_{\ti{g}} |\alpha|^2_{\ti{g}}\geq \frac{1}{2}|\ti{\alpha}_{i,j}|^2+|\ti{\alpha}_{i,\ov{j}}|^2+\frac{1}{2}\ti{g}(\Delta^H_{\ti{g}} a,a)-C\tr{g}{\ti{g}}|\alpha|^2_{\ti{g}}-C|\nabla\ti{g}|_{g,\ti{g}}|\alpha|^2_{\ti{g}}.$$
We compute
$$d\alpha=d(\ti{\alpha}_i\ti{\theta}^i)=\ti{\alpha}_{i,j}\ti{\theta}^j\wedge\ti{\theta}^i+\ti{\alpha}_{i,\ov{j}}\ov{\ti{\theta}^j}\wedge\ti{\theta}^i+\ti{\alpha}_i\ti{\Theta}^i,$$
and recall that
$$da=d\alpha+\ov{d\alpha}=\ti{\omega}-\omega,$$
and taking the $(1,1)$-part
$$(da)^{(1,1)}=\ti{\omega}-\omega^{(1,1)},$$
hence
$$|(da)^{(1,1)}|^2_{\ti{g}}=2-2\tr{\ti{g}}{g}+|g|^2_{\ti{g}}\geq 2-2\tr{\ti{g}}{g}+\frac{(\tr{\ti{g}}{g})^2}{2}\geq \frac{(\tr{\ti{g}}{g})^2}{4}-2,$$
but we also have
$$(da)^{(1,1)}=(\ti{\alpha}_{i,\ov{j}}-\ov{\ti{\alpha}_{j,\ov{i}}})\ov{\ti{\theta}^j}\wedge\ti{\theta}^i,$$
$$|(da)^{(1,1)}|^2_{\ti{g}}=|\ti{\alpha}_{i,\ov{j}}-\ov{\ti{\alpha}_{j,\ov{i}}}|^2\leq 4|\ti{\alpha}_{i,\ov{j}}|^2,$$
and so
$$|\ti{\alpha}_{i,\ov{j}}|^2\geq \frac{(\tr{\ti{g}}{g})^2}{16}-\frac{1}{2}\geq \frac{(\tr{g}{\ti{g}})^2}{C}-\frac{1}{2}.$$
This gives
$$\Delta_{\ti{g}} |\alpha|^2_{\ti{g}}\geq \frac{1}{2}|\ti{\alpha}_{i,j}|^2+\frac{1}{2}|\ti{\alpha}_{i,\ov{j}}|^2+\frac{(\tr{g}{\ti{g}})^2}{C}+\frac{1}{2}\ti{g}(\Delta^H_{\ti{g}} a,a)-C\tr{g}{\ti{g}}|\alpha|^2_{\ti{g}}-C|\nabla\ti{g}|_{g,\ti{g}}|\alpha|^2_{\ti{g}}-C,$$
and we also have
$$\frac{1}{2}|\ti{\alpha}_{i,j}|^2+\frac{1}{2}|\ti{\alpha}_{i,\ov{j}}|^2\geq\frac{1}{4}|\ti{\nabla}a|^2_{\ti{g}}+\frac{1}{4}|\ov{\ti{\nabla}}a|^2_{\ti{g}},$$
and recalling that $|a|^2_{\ti{g}}=2|\alpha|^2_{\ti{g}},$ we obtain
\begin{equation}\label{drei}
\Delta_{\ti{g}} |a|^2_{\ti{g}}\geq \frac{1}{2}|\ti{\nabla}a|^2_{\ti{g}}+\frac{1}{2}|\ov{\ti{\nabla}}a|^2_{\ti{g}}+\frac{(\tr{g}{\ti{g}})^2}{C}+\ti{g}(\Delta^H_{\ti{g}} a,a)-C\tr{g}{\ti{g}}|a|^2_{\ti{g}}-C|\nabla\ti{g}|_{g,\ti{g}}|a|^2_{\ti{g}}-C.
\end{equation}
We now deal with the term with the Hodge Laplacian. Using the gauge-fixing condition \eqref{gauge} we have
$$\Delta^H_{\ti{g}} a=-dd^*_{\ti{g}}a-d^*_{\ti{g}}da=-d^*_{\ti{g}}da=d^*_{\ti{g}}\omega-d^*_{\ti{g}}\ti{\omega}=d^*_{\ti{g}}\omega,$$
since $*_{\ti{g}}\ti{\omega}=\ti{\omega}$ and so $d^*_{\ti{g}}\ti{\omega}=0$.
We also have the well-known formula
$$*_{\ti{g}}\omega=\frac{2\omega\wedge\ti{\omega}}{\ti{\omega}^2}\ti{\omega}-\omega+2\omega^{(2,0)+(0,2)}=\tr{\ti{g}}{g}\,\ti{\omega}-\omega+2\omega^{(2,0)+(0,2)},$$
$$d*_{\ti{g}}\omega=d\tr{\ti{g}}{g}\wedge\ti{\omega}+2d(\omega^{(2,0)+(0,2)}),$$
$$*_{\ti{g}}d*_{\ti{g}}\omega=d^c\tr{\ti{g}}{g}+2*_{\ti{g}}d(\omega^{(2,0)+(0,2)}),$$
and so
$$\Delta^H_{\ti{g}} a=d^*_{\ti{g}}\omega=-d^c\tr{\ti{g}}{g}-2*_{\ti{g}}d(\omega^{(2,0)+(0,2)}),$$
\[\begin{split}
|\ti{g}(\Delta^H_{\ti{g}} a,a)|&\leq |\ti{g}(d^c\tr{\ti{g}}{g},a)|+2|\ti{g}(*_{\ti{g}}d(\omega^{(2,0)+(0,2)}),a)|\\
&\leq C|\nabla \tr{\ti{g}}{g}|_{\ti{g}} |a|_{\ti{g}}+2|*_{\ti{g}}d(\omega^{(2,0)+(0,2)})|_{\ti{g}}|a|_{\ti{g}}\\
&=C|\nabla (\tr{g}{\ti{g}}e^{-F})|_{\ti{g}} |a|_{\ti{g}}+2|d(\omega^{(2,0)+(0,2)})|_{\ti{g}}|a|_{\ti{g}}\\
&\leq C|\nabla \tr{g}{\ti{g}}|_{\ti{g}} |a|_{\ti{g}}+C\tr{g}{\ti{g}}|\nabla (e^{-F})|_{\ti{g}} |a|_{\ti{g}}+C(\tr{g}{\ti{g}})^{\frac{3}{2}} |a|_{\ti{g}}\\
&\leq C|\nabla \tr{g}{\ti{g}}|_{\ti{g}} |a|_{\ti{g}}+C(\tr{g}{\ti{g}})^{\frac{3}{2}} |a|_{\ti{g}}\\
&\leq \ve \frac{|\nabla \tr{g}{\ti{g}}|^2_{\ti{g}}}{\tr{g}{\ti{g}}}+\frac{C}{\ve} \tr{g}{\ti{g}}|a|_{\ti{g}}^2+\ve(\tr{g}{\ti{g}})^{2}+\frac{C}{\ve}|a|_{\ti{g}}^4\\
&\leq \ve \frac{|\nabla \tr{g}{\ti{g}}|^2_{\ti{g}}}{\tr{g}{\ti{g}}}+2\ve(\tr{g}{\ti{g}})^{2}+\frac{C}{\ve^2}|a|_{\ti{g}}^4\\
&\leq \ve |\nabla\ti{g}|^2_{g,\ti{g}}+2\ve(\tr{g}{\ti{g}})^{2}+\frac{C}{\ve^2}|a|_{\ti{g}}^4,
\end{split}\]
using \eqref{ein1}, and so assuming without loss that $\ve$ is small, and substituting in \eqref{drei} gives
$$\Delta_{\ti{g}} |a|^2_{\ti{g}}\geq \frac{1}{2}|\ti{\nabla}a|^2_{\ti{g}}+\frac{1}{2}|\ov{\ti{\nabla}}a|^2_{\ti{g}}+\frac{(\tr{g}{\ti{g}})^2}{C}-\ve |\nabla\ti{g}|^2_{g,\ti{g}}-\frac{C}{\ve^2}|a|_{\ti{g}}^4-C,$$
which is exactly \eqref{zwei}.\\

\noindent {\bf Proof of the main claim \eqref{vier} modulo \eqref{dolore1} and \eqref{dolore2}.}
We now need to prove the main claim \eqref{vier}. Since the Hodge Laplacian is a real operator, we have
$$\Delta^H_{\ti{g}} a=\Delta^H_{\ti{g}} \alpha+\Delta^H_{\ti{g}} \ov{\alpha}=\Delta^H_{\ti{g}} \alpha+\ov{\Delta^H_{\ti{g}} \alpha},$$
and
$$(\Delta^H_{\ti{g}} \ov{\alpha})^{(0,1)}=\ov{(\Delta^H_{\ti{g}} \alpha)^{(1,0)}},$$
and so
\begin{equation}\label{longa}\begin{split}
\ti{g}(\Delta^H_{\ti{g}} a,a)&=\ti{g}(\Delta^H_{\ti{g}} \alpha+\Delta^H_{\ti{g}} \ov{\alpha},\alpha+\ov{\alpha})\\
&=\ti{g}(\Delta^H_{\ti{g}} \alpha,\ov{\alpha})+\ti{g}(\Delta^H_{\ti{g}} \alpha,\alpha)+\ti{g}(\Delta^H_{\ti{g}} \ov{\alpha},\alpha)+\ti{g}(\Delta^H_{\ti{g}} \ov{\alpha},\ov{\alpha})\\
&=\ti{g}((\Delta^H_{\ti{g}} \alpha)^{(1,0)},\ov{\alpha})+\ti{g}((\Delta^H_{\ti{g}} \alpha)^{(0,1)},\alpha)+\ti{g}((\Delta^H_{\ti{g}} \ov{\alpha})^{(0,1)},\alpha)+\ti{g}((\Delta^H_{\ti{g}} \ov{\alpha})^{(1,0)},\ov{\alpha})\\
&=\ti{g}((\Delta^H_{\ti{g}} \alpha)^{(1,0)},\ov{\alpha})+\ti{g}((\Delta^H_{\ti{g}} \alpha)^{(0,1)},\alpha)+\ti{g}(\ov{(\Delta^H_{\ti{g}} \alpha)^{(1,0)}},\alpha)+\ti{g}(\ov{(\Delta^H_{\ti{g}} \alpha)^{(0,1)}},\ov{\alpha})\\
&=2\mathrm{Re}\ti{g}((\Delta^H_{\ti{g}} \alpha)^{(1,0)},\ov{\alpha})+2\mathrm{Re}\ti{g}((\Delta^H_{\ti{g}} \alpha)^{(0,1)},\alpha).
\end{split}\end{equation}
Next, we claim that
\begin{equation}\label{dolore1}
(\Delta^H_{\ti{g}} \alpha)^{(1,0)}=2\ti{\alpha}_{k,\ov{i}i}\ti{\theta}^k,
\end{equation}
\begin{equation}\label{dolore2}
(\Delta^H_{\ti{g}} \alpha)^{(0,1)}=\left(-4\ti{\alpha}_{i,j}\ti{N}^j_{\ov{k}\ov{i}}-2\ti{\alpha}_j \ti{K}^j_{i\ov{k}\ov{i}}+2\ti{\alpha}_i\ti{N}^i_{\ov{j}\ov{k},j}\right)\ov{\ti{\theta}^k}.
\end{equation}
We will prove these by long direct computations, but first let us assume we have these and complete the proof of the main claim \eqref{vier}. Combining \eqref{longa}, \eqref{dolore1} and \eqref{dolore2} gives
\begin{equation}\label{lap}\begin{split}
\frac{1}{2}\ti{g}(\Delta^H_{\ti{g}} a,a)&=2\mathrm{Re}(\ti{\alpha}_{k,\ov{i}i}\ov{\ti{\alpha}_k})+\mathrm{Re}\left(-4\ti{\alpha}_{i,j}\ti{\alpha}_k\ti{N}^j_{\ov{k}\ov{i}}-2\ti{\alpha}_j\ti{\alpha}_k\ti{K}^j_{i\ov{k}\ov{i}}+2\ti{\alpha}_i\ti{\alpha}_k\ti{N}^i_{\ov{j}\ov{k},j}\right),
\end{split}\end{equation}
and we can bound the last $3$ terms in \eqref{lap} using the same method as before. For the first term, using the same notation as earlier, we can write
$$\ti{\alpha}_{i,j}\ti{\alpha}_k\ti{N}^j_{\ov{k}\ov{i}}=\sum_{i,k=1}^2\lambda_i^{-1}\lambda_k^{-1}\alpha_{i,j}\alpha_kN^j_{\ov{k}\ov{i}}=\sum_{i\neq k}\lambda_i^{-1}\lambda_k^{-1}\alpha_{i,j}\alpha_kN^j_{\ov{k}\ov{i}}=e^{-F}\sum_{i\neq k}\alpha_{i,j}\alpha_kN^j_{\ov{k}\ov{i}},$$
using that $N^j_{\ov{k}\ov{i}}$ is skew-symmetric in $i,k$ as well as the Calabi-Yau equation \eqref{cy}. Using the Calabi-Yau equation again, we can bound
\[\begin{split}
-\mathrm{Re}(4\ti{\alpha}_{i,j}\ti{\alpha}_k\ti{N}^j_{\ov{k}\ov{i}})&\leq C\sum_j\sum_{i\neq k}|\alpha_{i,j}||\alpha_k|
\leq C\sum_j\sum_{i\neq k}\lambda_i^{-\frac{1}{2}}\lambda_k^{-\frac{1}{2}}|\alpha_{i,j}||\alpha_k|\lambda_j^{-\frac{1}{2}}\lambda_j^{\frac{1}{2}}\\
&=C\sum_j\sum_{i\neq k}|\ti{\alpha}_{i,j}||\ti{\alpha}_k|\lambda_j^{\frac{1}{2}}\leq C|\alpha|_{\ti{g}}\sqrt{\tr{g}{\ti{g}}}\sum_{i,j}|\ti{\alpha}_{i,j}|\\
&\leq \frac{1}{2}|\ti{\alpha}_{i,j}|^2+C\tr{g}{\ti{g}}|\alpha|_{\ti{g}}^2.
\end{split}\]
The second terms in the parenthesis in \eqref{lap} is bounded by $C\tr{g}{\ti{g}}|\alpha|^2_{\ti{g}}$ arguing exactly as we did to bound \eqref{list}: using the skew-symmetry of $K^j_{i\ov{k}\ov{i}}$ in $k,i$ and the Calabi-Yau equation \eqref{cy} we have
$$\ti{\alpha}_j\ti{\alpha}_k\ti{K}^j_{i\ov{k}\ov{i}}=\sum_{i,k=1}^2\lambda_i^{-1}\lambda_k^{-1}\alpha_j\alpha_kK^j_{i\ov{k}\ov{i}}=\sum_{i\neq k}\lambda_i^{-1}\lambda_k^{-1}\alpha_j\alpha_kK^j_{i\ov{k}\ov{i}}=
e^{-F}\sum_{i\neq k}\alpha_j\alpha_kK^j_{i\ov{k}\ov{i}},$$
and so
$$-\mathrm{Re}(2\ti{\alpha}_j\ti{\alpha}_k\ti{K}^j_{i\ov{k}\ov{i}})\leq C|\alpha|^2_g\leq C\tr{g}{\ti{g}}|\alpha|^2_{\ti{g}}.$$
The last term in the parenthesis in \eqref{lap}
\begin{equation}\label{ddd}
\ti{\alpha}_i\ti{\alpha}_k\ti{N}^i_{\ov{j}\ov{k},j},
\end{equation}
requires more work. First, we explain in detail the term $|\nabla\ti{g}|^2_{g,\ti{g}}$ in \eqref{ein}. As in \cite[(3.9)]{TWY} there are functions $a^i_{k\ell}$ defined by
$$da^i_m-a^i_j\theta^j_m+a^k_m\ti{\theta}^i_k=a^i_{k\ell}a^k_m\ti{\theta}^\ell,$$
which are the components of $\nabla\ti{g}$ in the frame $\{\ti{\theta}^i\}$, so that \cite[Lemma 4.2]{TWY}
$$|\nabla\ti{g}|^2_{\ti{g}}=|a^i_{k\ell}|^2.$$
The mixed norm $|\nabla\ti{g}|^2_{g,\ti{g}}$ that appears in \eqref{ein} is given by
\begin{equation}\label{mixed}
|\nabla\ti{g}|^2_{g,\ti{g}}=|a^i_{p\ell}a^p_k|^2.
\end{equation}
We can now go back to the term in \eqref{ddd} and using \cite[Lemma 4.4 (i)]{TWY} we can write it as
$$\ti{\alpha}_i\ti{\alpha}_k\ti{N}^i_{\ov{j}\ov{k},j}=\alpha_\ell \alpha_m b^m_k \ov{b^s_k} b^u_j \ov{b^r_j} N^\ell_{\ov{r}\ov{s},u}+
\alpha_\ell \alpha_m b^\ell_i b^m_k \ov{b^s_k} \ov{b^r_j} N^t_{\ov{r}\ov{s}} a^i_{uj}a^u_t,$$
and using the skew-symmetry of $N^i_{\ov{j}\ov{k},j}$ in $j,k$ and the Calabi-Yau equation we have
$$\alpha_\ell \alpha_m b^m_k \ov{b^s_k} b^u_j \ov{b^r_j} N^\ell_{\ov{r}\ov{s},u}=\sum_{j,k=1}^2\lambda_j^{-1}\lambda_k^{-1}\alpha_i\alpha_kN^i_{\ov{j}\ov{k},j}=\sum_{j\neq k}^2\lambda_j^{-1}\lambda_k^{-1}\alpha_i\alpha_kN^i_{\ov{j}\ov{k},j}
=e^{-F}\sum_{j\neq k}\alpha_i\alpha_kN^i_{\ov{j}\ov{k},j},$$
whose absolute value is bounded by $ C|\alpha|^2_g\leq C\tr{g}{\ti{g}}|\alpha|^2_{\ti{g}}.$ Similarly, using the skew-symmetry of $N^t_{\ov{j}\ov{k}}$ in $j,k$ and the Calabi-Yau equation we have
\[\begin{split}
\alpha_\ell \alpha_m b^\ell_i b^m_k \ov{b^s_k} \ov{b^r_j} N^t_{\ov{r}\ov{s}} a^i_{uj}a^u_t&=\sum_{i,j,k=1}^2\alpha_i \alpha_k \lambda_i^{-\frac{1}{2}}\lambda_k^{-1} \lambda_j^{-\frac{1}{2}} N^t_{\ov{j}\ov{k}} a^i_{uj}a^u_t\\
&=\sum_i \sum_{j\neq k}\alpha_i \alpha_k \lambda_i^{-\frac{1}{2}}\lambda_k^{-1} \lambda_j^{-\frac{1}{2}} N^t_{\ov{j}\ov{k}} a^i_{uj}a^u_t\\
&=e^{-\frac{F}{2}}\sum_i \sum_{j\neq k}\alpha_i \alpha_k \lambda_i^{-\frac{1}{2}}\lambda_k^{-\frac{1}{2}}  N^t_{\ov{j}\ov{k}} a^i_{uj}a^u_t\\
&=e^{-\frac{F}{2}}\sum_i \sum_{j\neq k}\ti{\alpha}_i \ti{\alpha}_k   N^t_{\ov{j}\ov{k}} a^i_{uj}a^u_t,
\end{split}\]
whose absolute value is bounded by
$$C|\alpha|^2_{\ti{g}}|a^i_{uj}a^u_t|=C|\alpha|^2_{\ti{g}}|\nabla\ti{g}|_{g,\ti{g}},$$
and so
$$-\mathrm{Re}(2\ti{\alpha}_i\ti{\alpha}_k\ti{N}^i_{\ov{j}\ov{k},j})\leq C\tr{g}{\ti{g}}|\alpha|^2_{\ti{g}}+C|\alpha|^2_{\ti{g}}|\nabla\ti{g}|_{g,\ti{g}}.$$
Putting all these together we obtain
$$\frac{1}{2}\ti{g}(\Delta^H_{\ti{g}} a,a)\leq 2\mathrm{Re}(\ti{\alpha}_{k,\ov{i}i}\ov{\ti{\alpha}_k})+C\tr{g}{\ti{g}}|\alpha|^2_{\ti{g}}+C|\alpha|^2_{\ti{g}}|\nabla\ti{g}|_{g,\ti{g}}+\frac{1}{2}|\ti{\alpha}_{i,j}|^2,$$
which is exactly the main claim \eqref{vier}.\\

\noindent {\bf Proof of \eqref{dolore1} and \eqref{dolore2}.}
To complete the proof of the main claim \eqref{vier}, we are left with showing \eqref{dolore1} and \eqref{dolore2}, which are Bochner-Kodaira-Weitzenb\"ock type formulas. For these, recall that by definition of the Hodge star $*_{\ti{g}}$, for any two $(p,q)$-forms $\beta,\gamma$ we have
$$\beta\wedge \ov{*_{\ti{g}}\gamma}=\ti{g}(\beta,\gamma)\ti{\omega}^2=\ti{g}(\beta,\gamma)\sqrt{-1}\ti{\theta}^1\wedge\ov{\ti{\theta}^1}\wedge\sqrt{-1}\ti{\theta}^2\wedge\ov{\ti{\theta}^2}
=-\ti{g}(\beta,\gamma)\ti{\theta}^1\wedge\ov{\ti{\theta}^1}\wedge\ti{\theta}^2\wedge\ov{\ti{\theta}^2},$$
and also $*_{\ti{g}}\ov{\alpha}=\ov{*_{\ti{g}}\alpha}$ (i.e. $*_{\ti{g}}$ is a real operator),
from which we can compute the Hodge star on basic combinations of our frame elements, using the following notation: for $i\in \{1,2\}$ we let $\hat{\imath}\in\{1,2\}$ be such that $\{i,\hat{\imath}\}=\{1,2\}$ (unordered), i.e. $\hat{1}=2,\hat{2}=1$. Then from the definition we have
$$*_{\ti{g}}(\ti{\theta}^i\wedge\ov{\ti{\theta}^i}\wedge\ti{\theta}^{\hat{\imath}}\wedge\ov{\ti{\theta}^{\hat{\imath}}})=-1,$$
$$*_{\ti{g}}\ti{\theta}^i=\ti{\theta}^i\wedge\ti{\theta}^{\hat{\imath}}\wedge\ov{\ti{\theta}^{\hat{\imath}}},\quad *_{\ti{g}}\ov{\ti{\theta}^i}=-\ov{\ti{\theta}^i}\wedge\ti{\theta}^{\hat{\imath}}\wedge\ov{\ti{\theta}^{\hat{\imath}}}$$
$$*_{\ti{g}}(\ti{\theta}^i\wedge \ti{\theta}^{\hat{\imath}}\wedge\ov{ \ti{\theta}^{\hat{\imath}}})=-\ti{\theta}^i,\quad *_{\ti{g}}(\ti{\theta}^i\wedge \ti{\theta}^{\hat{\imath}}\wedge\ov{ \ti{\theta}^{i}})=\ti{\theta}^{\hat{\imath}},$$
$$*_{\ti{g}}(\ti{\theta}^{\hat{\imath}}\wedge \ov{\ti{\theta}^{\hat{\imath}}}\wedge\ov{ \ti{\theta}^{i}})=\ov{\ti{\theta}^{i}},\quad *_{\ti{g}}(\ti{\theta}^i\wedge \ov{\ti{\theta}^{\hat{\imath}}}\wedge\ov{ \ti{\theta}^{i}})=-\ov{\ti{\theta}^{\hat{\imath}}},$$
$$*_{\ti{g}}(\ti{\theta}^i\wedge\ti{\theta}^j)=\ti{\theta}^i\wedge\ti{\theta}^j,\quad *_{\ti{g}}(\ov{\ti{\theta}^i}\wedge\ov{\ti{\theta}^j})=\ov{\ti{\theta}^i}\wedge\ov{\ti{\theta}^j},$$
$$*_{\ti{g}}(\ti{\theta}^i\wedge\ov{\ti{\theta}^j})=\begin{cases}\ti{\theta}^{\hat{\imath}}\wedge \ov{\ti{\theta}^{\hat{\imath}}}\quad\quad\text{if }i=j,\\
-\ti{\theta}^i\wedge\ov{\ti{\theta}^j}\quad\,\text{if }i=\hat{\jmath}.\end{cases}$$
With these, we can now start the computation of $\Delta^H_{\ti{g}}\alpha=-dd^*_{\ti{g}}\alpha-d^*_{\ti{g}}d\alpha$, recalling that
$$d^*_{\ti{g}}\alpha=-*_{\ti{g}}d*_{\ti{g}}\alpha,$$
we compute
$$*_{\ti{g}}\alpha=\ti{\alpha}_i\ti{\theta}^i\wedge\ti{\theta}^{\hat{\imath}}\wedge\ov{\ti{\theta}^{\hat{\imath}}},$$
$$d*_{\ti{g}}\alpha=\ti{\alpha}_{i,\ov{j}}\ov{\ti{\theta}^j}\wedge \ti{\theta}^i\wedge\ti{\theta}^{\hat{\imath}}\wedge\ov{\ti{\theta}^{\hat{\imath}}}=\ti{\alpha}_{i,\ov{i}}\ov{\ti{\theta}^i}\wedge \ti{\theta}^i\wedge\ti{\theta}^{\hat{\imath}}\wedge\ov{\ti{\theta}^{\hat{\imath}}},$$
$$d^*_{\ti{g}}\alpha=-\ti{\alpha}_{i,\ov{i}},$$
\begin{equation}\label{first}
dd^*_{\ti{g}}\alpha=-\ti{\alpha}_{i,\ov{i}k}\ti{\theta}^k-\ti{\alpha}_{i,\ov{i}\ov{k}}\ov{\ti{\theta}^k}.
\end{equation}
$$d\alpha=\ti{\alpha}_{i,j}\ti{\theta}^j\wedge\ti{\theta}^i+\ti{\alpha}_{i,\ov{j}}\ov{\ti{\theta}^j}\wedge\ti{\theta}^i+\ti{\alpha}_i\ti{N}^i_{\ov{j}\ov{k}}\ov{\ti{\theta}^j}\wedge\ov{\ti{\theta}^k}$$
\begin{equation}\label{last}
*_{\ti{g}}d\alpha=\ti{\alpha}_{i,j}\ti{\theta}^j\wedge\ti{\theta}^i-\ti{\alpha}_{i,\ov{i}}\ti{\theta}^{\hat{\imath}}\wedge\ov{\ti{\theta}^{\hat{\imath}}}
+\ti{\alpha}_{i,\ov{\hat{\imath}}}\ti{\theta}^{i}\wedge\ov{\ti{\theta}^{\hat{\imath}}}+\ti{\alpha}_i\ti{N}^i_{\ov{j}\ov{k}}\ov{\ti{\theta}^j}\wedge\ov{\ti{\theta}^k}.
\end{equation}
We first prove \eqref{dolore1}. We are thus interested in the $(1,0)$-part of $d^*_{\ti{g}}d\alpha$, so we first take the exterior derivative of \eqref{last} and take its $(2,1)$-part (which becomes $(1,0)$ after taking the Hodge star)
\begin{equation}\label{lost}
(d*_{\ti{g}}d\alpha)^{(2,1)}=\ti{\alpha}_{i,j\ov{k}}\ov{\ti{\theta}^k}\wedge\ti{\theta}^j\wedge\ti{\theta}^i-\ti{\alpha}_{i,\ov{i}k}\ti{\theta}^k\wedge\ti{\theta}^{\hat{\imath}}\wedge\ov{\ti{\theta}^{\hat{\imath}}}
+\ti{\alpha}_{i,\ov{\hat{\imath}}k}\ti{\theta}^k\wedge\ti{\theta}^{i}\wedge\ov{\ti{\theta}^{\hat{\imath}}}+2\ti{\alpha}_i\ti{N}^i_{\ov{j}\ov{k}}\ov{\ti{N}^j_{\ov{p}\ov{q}}}\ti{\theta}^p\wedge\ti{\theta}^q\wedge\ov{\ti{\theta}^k},
\end{equation}
and in the first term in \eqref{lost} we must have $j=\hat{\imath}$ and we can write it as
$$\ti{\alpha}_{i,{\hat{\imath}}\ov{i}}\ov{\ti{\theta}^i}\wedge\ti{\theta}^{\hat{\imath}}\wedge\ti{\theta}^i+\ti{\alpha}_{i,{\hat{\imath}}\ov{\hat{\imath}}}\ov{\ti{\theta}^{\hat{\imath}}}\wedge\ti{\theta}^{\hat{\imath}}\wedge\ti{\theta}^i,$$
and its Hodge star equals
$$-\ti{\alpha}_{i,{\hat{\imath}}\ov{i}}\ti{\theta}^{\hat{\imath}}+\ti{\alpha}_{i,\hat{\imath}\ov{\hat{\imath}}}\ti{\theta}^i=
(-\ti{\alpha}_{\hat{k},k\ov{\hat{k}}}+\ti{\alpha}_{k,\hat{k}\ov{\hat{k}}})\ti{\theta}^k=(-\ti{\alpha}_{i,k\ov{i}}+\ti{\alpha}_{k,i\ov{i}})\ti{\theta}^k.$$
In the second term we must have $k=i$ and in the third term $k=\hat{\imath}$ so we can write them as
$$-\ti{\alpha}_{i,\ov{i}i}\ti{\theta}^i\wedge\ti{\theta}^{\hat{\imath}}\wedge\ov{\ti{\theta}^{\hat{\imath}}}
+\ti{\alpha}_{i,\ov{\hat{\imath}}\hat{\imath}}\ti{\theta}^{\hat{\imath}}\wedge\ti{\theta}^{i}\wedge\ov{\ti{\theta}^{\hat{\imath}}},$$
and their Hodge star equals
$$(\ti{\alpha}_{i,\ov{i}i}+\ti{\alpha}_{i,\ov{\hat{\imath}}\hat{\imath}})\ti{\theta}^i=\ti{\alpha}_{k,\ov{i}i}\ti{\theta}^k.$$
In the last term in \eqref{lost} we must have $q=\hat{p}$ so it equals
$$2\ti{\alpha}_i\ti{N}^i_{\ov{j}\ov{p}}\ov{\ti{N}^j_{\ov{p}\ov{\hat{p}}}}\ti{\theta}^p\wedge\ti{\theta}^{\hat{p}}\wedge\ov{\ti{\theta}^p}
+2\ti{\alpha}_i\ti{N}^i_{\ov{j}\ov{\hat{p}}}\ov{\ti{N}^j_{\ov{p}\ov{\hat{p}}}}\ti{\theta}^p\wedge\ti{\theta}^{\hat{p}}\wedge\ov{\ti{\theta}^{\hat{p}}},$$
and its Hodge star equals
\[\begin{split}
2\ti{\alpha}_i\ti{N}^i_{\ov{j}\ov{p}}\ov{\ti{N}^j_{\ov{p}\ov{\hat{p}}}}\ti{\theta}^{\hat{p}}
-2\ti{\alpha}_i\ti{N}^i_{\ov{j}\ov{\hat{p}}}\ov{\ti{N}^j_{\ov{p}\ov{\hat{p}}}}\ti{\theta}^p&=\left(2\ti{\alpha}_i\ti{N}^i_{\ov{j}\ov{\hat{k}}}\ov{\ti{N}^j_{\ov{\hat{k}}\ov{k}}}-
2\ti{\alpha}_i\ti{N}^i_{\ov{j}\ov{\hat{k}}}\ov{\ti{N}^j_{\ov{k}\ov{\hat{k}}}}\right)\ti{\theta}^k\\
&=\left(2\ti{\alpha}_i\ti{N}^i_{\ov{j}\ov{\ell}}\ov{\ti{N}^j_{\ov{\ell}\ov{k}}}-
2\ti{\alpha}_i\ti{N}^i_{\ov{j}\ov{\ell}}\ov{\ti{N}^j_{\ov{k}\ov{\ell}}}\right)\ti{\theta}^k\\
&=-4\ti{\alpha}_i\ti{N}^i_{\ov{j}\ov{\ell}}\ov{\ti{N}^j_{\ov{k}\ov{\ell}}}\ti{\theta}^k,\end{split}\]
using the skew-symmetry of $\ti{N}^j_{\ov{k}\ov{\ell}}$ in $k,\ell$. Putting these together gives
\[\begin{split}
-(d^*_{\ti{g}}d\alpha)^{(1,0)}&=(*_{\ti{g}}d*_{\ti{g}}d\alpha)^{(1,0)}=*_{\ti{g}}\left((d*_{\ti{g}}d\alpha)^{(2,1)}\right)\\
&=\left(-\ti{\alpha}_{i,k\ov{i}}+\ti{\alpha}_{k,i\ov{i}}+\ti{\alpha}_{k,\ov{i}i}-4\ti{\alpha}_i\ti{N}^i_{\ov{j}\ov{\ell}}\ov{\ti{N}^j_{\ov{k}\ov{\ell}}}\right)\ti{\theta}^k,
\end{split}\]
and combining this with \eqref{first} gives
$$(\Delta^H_{\ti{g}}\alpha)^{(1,0)}=\left(\ti{\alpha}_{i,\ov{i}k}-\ti{\alpha}_{i,k\ov{i}}+\ti{\alpha}_{k,i\ov{i}}+\ti{\alpha}_{k,\ov{i}i}-4\ti{\alpha}_i\ti{N}^i_{\ov{j}\ov{\ell}}\ov{\ti{N}^j_{\ov{k}\ov{\ell}}}\right)\ti{\theta}^k.$$
From the commutation relation \eqref{comm1} we obtain
$$\ti{\alpha}_{i,\ov{i}k}-\ti{\alpha}_{i,k\ov{i}}=-\ti{\alpha}_j\ti{R}^j_{ik\ov{i}},$$
$$\ti{\alpha}_{k,i\ov{i}}-\ti{\alpha}_{k,\ov{i}i}=\ti{\alpha}_j\ti{R}^j_{ki\ov{i}},$$
while \cite[(2.16)]{TWY} gives
$$\ti{R}^j_{ki\ov{i}}-\ti{R}^j_{ik\ov{i}}=4\ti{N}^j_{\ov{p}\ov{i}}\ov{\ti{N}^p_{\ov{k}\ov{i}}},$$
and so
$$(\Delta^H_{\ti{g}}\alpha)^{(1,0)}=\left(2\ti{\alpha}_{k,\ov{i}i}+4\ti{\alpha}_j\ti{N}^j_{\ov{p}\ov{i}}\ov{\ti{N}^p_{\ov{k}\ov{i}}}-4\ti{\alpha}_i\ti{N}^i_{\ov{j}\ov{\ell}}\ov{\ti{N}^j_{\ov{k}\ov{\ell}}}\right)\ti{\theta}^k
=2\ti{\alpha}_{k,\ov{i}i}\ti{\theta}^k,$$
which proves \eqref{dolore1}. The proof of \eqref{dolore2} is similar. First, from \eqref{first} we have
\begin{equation}\label{first1}
(dd^*_{\ti{g}}\alpha)^{(0,1)}=-\ti{\alpha}_{i,\ov{i}\ov{k}}\ov{\ti{\theta}^k}.
\end{equation}
We are then interested in the $(0,1)$-part of $d^*_{\ti{g}}d\alpha$, so we take the exterior derivative of \eqref{last} and then take the $(1,2)$-part to obtain
\begin{equation}\label{lost2}\begin{split}
(d*_{\ti{g}}d\alpha)^{(1,2)}&=\ti{\alpha}_{i,j}\ti{N}^j_{\ov{p}\ov{q}}\ov{\ti{\theta}^p}\wedge\ov{\ti{\theta}^q}\wedge\ti{\theta}^i-\ti{\alpha}_{i,j}\ti{N}^i_{\ov{p}\ov{q}}\ti{\theta}^j\wedge\ov{\ti{\theta}^p}\wedge\ov{\ti{\theta}^q}
-\ti{\alpha}_{i,\ov{i}\ov{k}}\ov{\ti{\theta}^k}\wedge\ti{\theta}^{\hat{\imath}}\wedge\ov{\ti{\theta}^{\hat{\imath}}}\\
&+\ti{\alpha}_{i,\ov{\hat{\imath}}\ov{k}}\ov{\ti{\theta}^k}\wedge\ti{\theta}^{i}\wedge\ov{\ti{\theta}^{\hat{\imath}}}+\ti{\alpha}_{i,p}\ti{N}^i_{\ov{j}\ov{k}}\ti{\theta}^p\wedge\ov{\ti{\theta}^j}\wedge\ov{\ti{\theta}^k}
+\ti{\alpha}_i\ti{N}^i_{\ov{j}\ov{k},p}\ti{\theta}^p\wedge\ov{\ti{\theta}^j}\wedge\ov{\ti{\theta}^k}\\
&=\ti{\alpha}_{i,j}\ti{N}^j_{\ov{p}\ov{q}}\ov{\ti{\theta}^p}\wedge\ov{\ti{\theta}^q}\wedge\ti{\theta}^i
-\ti{\alpha}_{i,\ov{i}\ov{k}}\ov{\ti{\theta}^k}\wedge\ti{\theta}^{\hat{\imath}}\wedge\ov{\ti{\theta}^{\hat{\imath}}}
+\ti{\alpha}_{i,\ov{\hat{\imath}}\ov{k}}\ov{\ti{\theta}^k}\wedge\ti{\theta}^{i}\wedge\ov{\ti{\theta}^{\hat{\imath}}}\\
&+\ti{\alpha}_i\ti{N}^i_{\ov{j}\ov{k},p}\ti{\theta}^p\wedge\ov{\ti{\theta}^j}\wedge\ov{\ti{\theta}^k},
\end{split}\end{equation}
and in the first term in \eqref{lost2} we must have $q=\hat{p}$ so we can write it as
$$\ti{\alpha}_{p,j}\ti{N}^j_{\ov{p}\ov{\hat{p}}}\ov{\ti{\theta}^p}\wedge\ov{\ti{\theta}^{\hat{p}}}\wedge\ti{\theta}^p+\ti{\alpha}_{\hat{p},j}\ti{N}^j_{\ov{p}\ov{\hat{p}}}\ov{\ti{\theta}^p}\wedge\ov{\ti{\theta}^{\hat{p}}}\wedge\ti{\theta}^{\hat{p}},$$
and its Hodge star equals
$$\ti{\alpha}_{p,j}\ti{N}^j_{\ov{p}\ov{\hat{p}}}\ov{\ti{\theta}^{\hat{p}}}-\ti{\alpha}_{\hat{p},j}\ti{N}^j_{\ov{p}\ov{\hat{p}}}\ov{\ti{\theta}^p}=
\ti{\alpha}_{p,j}\ti{N}^j_{\ov{p}\ov{i}}\ov{\ti{\theta}^{i}}-\ti{\alpha}_{i,j}\ti{N}^j_{\ov{p}\ov{i}}\ov{\ti{\theta}^p}=-2\ti{\alpha}_{i,j}\ti{N}^j_{\ov{p}\ov{i}}\ov{\ti{\theta}^p}.$$
In the second and third terms in \eqref{lost2} we must have $k=i$ so we can write them as
$$-\ti{\alpha}_{i,\ov{i}\ov{i}}\ov{\ti{\theta}^i}\wedge\ti{\theta}^{\hat{\imath}}\wedge\ov{\ti{\theta}^{\hat{\imath}}}
+\ti{\alpha}_{i,\ov{\hat{\imath}}\ov{i}}\ov{\ti{\theta}^i}\wedge\ti{\theta}^{i}\wedge\ov{\ti{\theta}^{\hat{\imath}}},$$
and their Hodge star equals
$$-\ti{\alpha}_{i,\ov{i}\ov{i}}\ov{\ti{\theta}^i}-\ti{\alpha}_{i,\ov{\hat{\imath}}\ov{i}}\ov{\ti{\theta}^{\hat{\imath}}}=-\ti{\alpha}_{i,\ov{k}\ov{i}}\ov{\ti{\theta}^k}.$$
In the last term in \eqref{lost2} we must have $k=\hat{\jmath}$ and so we can write it as
$$\ti{\alpha}_i\ti{N}^i_{\ov{j}\ov{\hat{\jmath}},j}\ti{\theta}^j\wedge\ov{\ti{\theta}^j}\wedge\ov{\ti{\theta}^{\hat{\jmath}}}+\ti{\alpha}_i\ti{N}^i_{\ov{j}\ov{\hat{\jmath}},\hat{\jmath}}\ti{\theta}^{\hat{\jmath}}\wedge\ov{\ti{\theta}^j}\wedge\ov{\ti{\theta}^{\hat{\jmath}}},$$
and its Hodge star equals
$$\ti{\alpha}_i\ti{N}^i_{\ov{j}\ov{\hat{\jmath}},j}\ov{\ti{\theta}^{\hat{\jmath}}}-\ti{\alpha}_i\ti{N}^i_{\ov{j}\ov{\hat{\jmath}},\hat{\jmath}}\ov{\ti{\theta}^j}=
\ti{\alpha}_i\ti{N}^i_{\ov{j}\ov{k},j}\ov{\ti{\theta}^{k}}-\ti{\alpha}_i\ti{N}^i_{\ov{k}\ov{j},j}\ov{\ti{\theta}^k}=2\ti{\alpha}_i\ti{N}^i_{\ov{j}\ov{k},j}\ov{\ti{\theta}^{k}}.$$
Putting these together gives
\[\begin{split}-(d^*_{\ti{g}}d\alpha)^{(0,1)}&=(*_{\ti{g}}d*_{\ti{g}}d\alpha)^{(0,1)}=*_{\ti{g}}\left((d*_{\ti{g}}d\alpha)^{(1,2)}\right)\\
&=\left(-2\ti{\alpha}_{i,j}\ti{N}^j_{\ov{k}\ov{i}}-\ti{\alpha}_{i,\ov{k}\ov{i}}+2\ti{\alpha}_i\ti{N}^i_{\ov{j}\ov{k},j}\right)\ov{\ti{\theta}^{k}},
\end{split}\]
and combining this with \eqref{first1} gives
$$(\Delta^H_{\ti{g}}\alpha)^{(0,1)}=\left(\ti{\alpha}_{i,\ov{i}\ov{k}}-2\ti{\alpha}_{i,j}\ti{N}^j_{\ov{k}\ov{i}}-\ti{\alpha}_{i,\ov{k}\ov{i}}+2\ti{\alpha}_i\ti{N}^i_{\ov{j}\ov{k},j}\right)\ov{\ti{\theta}^{k}}.$$
From the commutation relation \eqref{comm2} we obtain
$$\ti{\alpha}_{i,\ov{i}\ov{k}}-\ti{\alpha}_{i,\ov{k}\ov{i}}=-2\ti{\alpha}_j\ti{K}^j_{i\ov{k}\ov{i}}-2\ti{\alpha}_{i,j}\ti{N}^j_{\ov{k}\ov{i}},$$
and so
$$(\Delta^H_{\ti{g}}\alpha)^{(0,1)}=\left(-4\ti{\alpha}_{i,j}\ti{N}^j_{\ov{k}\ov{i}}-2\ti{\alpha}_j\ti{K}^j_{i\ov{k}\ov{i}}+2\ti{\alpha}_i\ti{N}^i_{\ov{j}\ov{k},j}\right)\ov{\ti{\theta}^{k}},$$
which proves \eqref{dolore2} and concludes the proof of Proposition \ref{prop}, and hence of Theorem \ref{main}.

\end{document}